\documentstyle [12pt]{article}
\begin{document}

\title{ Kuhn-Tucker conditions for a convex programming problem in Banach spaces partially ordered by cone with empty interior.}
\author{ Feyzullah Ahmeto\u{g}lu\\
 Faculty of Education,
Giresun University,\\
 Giresun, Turkey. e-mail: feyzullah.ahmetoglu@giresun.edu.tr\\}
\date{}
\maketitle

{\bf 1. Introduction}\\

 During recent decades the theory of mathematical
programming in infinite dimensional spaces has been studied
extensively [1]-[7].

In order to obtain Kuhn-Tucker condition in mathematical
programming, problems usually are formulated in spaces where a
cone defining partial order has a nonempty interior. In these spaces  the existence of a
saddle point of the Lagrange function or Kuhn-Tucker conditions are established by
using of
some natural conditions like Slayter, regularity, etc. These well known methods fail in the cases
when the cones defining
partial order in the space have no interior points. 
$L_{p}[0,T]$ and $l_p (1 < p < \infty )$ spaces constitute examples for these cases.  In the present paper 
we explore spaces not necessarily having nonempty interior of the cone defining partial order. We  obtain a
differential form of Kuhn-Tucker conditions for a convex programming problem in Banach spaces without strong restriction
assuming the existence of nonempty interior of the cone defining partial order  in the space.\\

{\bf 2. Formulation of results.}\\

Let $X$ and $Y$ be reflexive Banach spaces partially ordered by
convex closed cones $K$ and $P$, respectively. A linear bounded
operator mapping $X$ into $Y$ we denote by $A$.

We investigate the problem of  minimization of the continuously
differentiable convex functional $I(x)$ under following additional
constraints:

$$ A x \leq b   \: \:  \: \: \: (b - A x) \in P$$

$$x \geq 0     \: \:  \: \: \: (x \in K)$$

The problem can be shortly formulated as

\begin{equation}
I(x) \rightarrow min
\end{equation}

\begin{equation}
A x \leq b \: \: \: \: \: x \geq 0
\end{equation}

{\bf Definition 1 }. We say that constraints (2) satisfy the
strong simultaneity condition, if there exists $\epsilon_0 > $
such that for every $\bar{b} \in \{ \bar{b}: || \: \bar{b} - b||
\leq \epsilon_0
\}$ the system $A x \leq \bar{b}, x \geq 0$ has a solution.\\

A point $p \in M$ is called an internal point of $M$, if for each
$z \in Y$ there exists a real number $\epsilon > 0$ such that for
each $\lambda $ satisfying $| \lambda | \leq \epsilon$ we have $p
+ \lambda z \in M$.\\

{\bf Lemma 1}. Suppose that the constraints (2) satisfy the strong
simultaneity condition. Then the set

$$M = \{ z \in Y : b - A x
\geq z, x \geq 0 \} $$

has internal points.\\

{\bf Proof}. In order to prove the lemma  it suffices to show that
a zero point is an internal point of $M$. In other words, for each
point $z \in Y, z \not= 0$, there exists a real number $\lambda',$
such that the constraints

$$b - A x \geq \lambda z, \: \: \: x \geq 0$$

are consistent for all $\lambda \in (0, \lambda')$. We choose
$\lambda' = \frac{\epsilon_0}{|| z ||}$. Then for each $\lambda
\in  (0, \lambda')$ we have $\lambda || z || < \epsilon_0$. Since
the conditions (2) are strongly  simultaneous, $b - A x \geq
\lambda z, x \geq 0$ for each  $\lambda \in (0, \lambda')$. The
proof is completed.\\

 {\bf Lemma 2}. Suppose that the constraints (2)
satisfy the strong simultaneity condition. Then the set

$$S= \{ (z,p) \in Y \times R: b - A x \geq z, I(x) \leq \rho, x
\geq 0 \}$$

has internal points. \\

{\bf Proof}. Clearly, there exists $x_0$ such that

\begin{equation}
b - Ax_0 \geq 0
\end{equation}

We show that $(0,I(x) +1)$ is an internal point of $S$. Let
$\rho_0 = I(x^{0})+1$. Let us show that for each $(z,\rho) \in Y
\times R$, there exists $\bar{\lambda}$, such that for arbitrary
$\lambda \in (0,\bar{\lambda})$ we have
 $ (\lambda z, \rho_0 + \lambda
\rho) \in S$. In other words, for arbitrary  $\lambda \in
(0,\bar{\lambda})$ there exists $x_{\lambda} \geq 0$, such that $b
- Ax_{\lambda} \geq \lambda z$ and $ \rho_0 + \lambda \rho \geq
I(x_{\lambda})$.

By Lemma 1 $0 \in Z$ is an internal point of $M$. Therefore, there
exist a real number $\lambda_0 > 0$ and a point $\bar{x}^{0}$ such
that

\begin{equation}
b - A\bar{x}_0 \geq \lambda_0 z
\end{equation}

By multiplying both sides of (3) by $1 -
\frac{\lambda}{\lambda_0}$ and both sides of (4) by
$\frac{\lambda}{\lambda_0}$ and  taking their sum, we get

$$b - A( \frac{\lambda}{\lambda_0}\bar{x}_0 + (1 - \frac{\lambda}{\lambda_0}
)x_0)) \geq \lambda z$$

Let  $x_{\lambda} = \frac{\lambda}{\lambda_0}\bar{x}_0 + (1 -
\frac{\lambda}{\lambda_0})x_0$. Then

\begin{equation}
b - A x_{\lambda} \geq \lambda z
\end{equation}

Since $I(x)$ is a convex functional we get

$$I(x_{\lambda}) \leq (1 - \frac{\lambda}{\lambda_0}) I(x^0) +
  \frac{\lambda}{\lambda_0}I(\bar{x}^0) = I(x_0) + \frac{\lambda}{\lambda_0}
  (I(\bar{x}_0) -I(x^0))$$

In order to  prove  $I(x_{\lambda}) \leq \rho_0 + \lambda \rho$ it
is enough to establish the following inequality

$$I(x^0) + \frac{\lambda}{\lambda_0}(I(\bar{x}_0) - I(x^0)) \leq
I(x^0) + 1 + \lambda \rho $$

The last inequality is held for all $\lambda \beta \leq 1$, where
$\beta = |\frac{I(\bar{x}_0 - I(x^0))}{\lambda_0} - \rho|$. Thus,
we can complete the proof by choosing $\bar{\lambda}$\\

$$\bar{\lambda}= \left \{
\begin{array}{ll} \displaystyle{min\{ \lambda_0, 1/\beta \}}
  & if \: \beta \not= 0\: \\
   \lambda_0  &  if \: \beta = 0\: \end{array}
 \right.$$ \\

 Lemma 2 is proved.\\

Let $X^{*}$ and $Y^{*}$ be the conjugate spaces of $X$ and $Y$,
respectively. The conjugate cone of $K$ is $K^{*}$:

$$K^{*} = \{x^{*} \in X^{*}: (x^{*},x) \geq 0 \:for \: all\: x \in K\}$$

The conjugate cone of P is defined similarly.

Let $X^{*}$ and $Y^{*}$ be partially ordered by $K^{*}$ and
$P^{*}$, respectively.\\

 {\bf Lemma 3}.
Suppose that the constraints (2) satisfy the strong simultaneity
condition. Then for any $z^{*} \in P^{*}, z^{*} \not= 0$, there
exists a point $x_{z^{*}} \geq 0$ such that

$$(z^{*}, b-Ax_{z^{*}})  > 0 $$

{\bf Proof}. For strong simultaneity of (2) for each $\xi \in Y,
||\xi|| \leq 1,$ there exists a point $x_{\xi} \geq 0$ such that

 $$b - A x_{\xi} \geq \epsilon_0 \xi$$

Let $z^{*} \in P^{*}$ and $z^{*} \not= 0$. Obviously, there exists
$z_0 \in Y$ such that

$$sup_{||z|| \leq 1} (z^{*},z) = (z^{*},z_0)$$

Since $||z_0|| \leq 1$, there exists $x_{z^{*}} \geq 0$ satisfying

$$b-Ax_{z^{*}}) \geq \epsilon_0 z_0$$

Now

$$(z^{*}, b - Ax_{z^{*}}) \geq \epsilon_0 (z^{*},z_0) = \epsilon_0
sup_{||z|| \leq 1}(z^{*},z) = \epsilon_0 ||z^{*}|| > 0$$

The lemma is proved.\\

The functional $L(x,z^{*}) = I(x) + (z^{*},Ax-b)$ is called a
Lagrange function.\\

 {\bf Definition 2 }. A pair $<x_0,z_0^{*}>$ is said to be a
saddle point of Lagrange function if $x_0 \geq 0, z_0^{*} \geq 0$
and for each $x \geq 0, z^{*} \geq 0$

\begin{equation}
L(x_0,z^{*}) \leq L(x_0,z_0^{*}) \leq L(x,z_0^{*})
\end{equation}

It can be easily shown that the existence of a saddle point of
Lagrange function implies the existence of a solution of problem
(1),(2). The inverse of this statement is also true:\\

{\bf Theorem  1}. Suppose that the constraints (2) satisfy the
strong simultaneity condition and the problem (1),(2) has a
solution $x_0$. Then there exists a non-zero linear functional
$z_0^{*}$ such that the pair $<x_0,z_0^{*}>$ is a saddle point of
Lagrange function.\\

{\bf Proof}. By Lemma 2 the set

$$S= \{ (z,\rho) \in Y \times R: b - A x \geq z, I(x) \leq \rho, x
\geq 0 \}$$

has internal points. By Lemma 3 for each $z^{*} \in P^{*}, z^{*}
\not= 0$, there exists a point $x_{z^{*}}$ such that
$(z^{*},b-Ax_{z^{*}}) > 0$. Thus, the strong simultaneity
condition implies both conditions of Theorem 1 of [1], which
states the existence of a saddle point.

Let us prove the existence of a saddle point in our case. Consider
the following sets in $Y \times R$

$$N = \{(z,\rho) \in Y \times R: z \geq 0, \rho \leq I(x_0)\}$$

$$N_1 = \{(z,\rho) \in Y \times R: z \geq 0, \rho < I(x_0)\}$$

The sets $S,N$ and $N_1$ are convex sets. Let us show that $S \cap
N_1 = \emptyset$. Indeed, if $x \geq 0$ and $Ax \leq b$, then for
all $(z,\rho) \in S$ we have $\rho \geq I(x) \geq I(x_0)$. On the
other hand, in $N_1$ $\rho < I(x_0)$. If $x \geq 0$ and $ b - Ax
\not\in P$, then in $N_1$ $z \geq 0$ but in S it is not held.
Done.

By Lemma 2, $S$ has an internal point. As a result, $S$ and $N_1$
are  disjoint convex sets and $S$ has an internal point.
Therefore, by well-known separation theorem [2], there exist
$(y_0^{*},\rho_0) \in Y^{*} \times R, (y_0^{*},\rho_0) \not= 0$
such that

\begin{equation}
\rho_0 \rho + (y_0^{*},z) \geq \rho_0 r + (y_0^{*},y)
\end{equation}

for all $(z,\rho) \in S$ and $(y,r)\in N_1$.

Since the closure of $N_1$ is $N$, (7) is also held for all $(y,r)
\in N$, which implies that $\rho_0 \geq 0$. Indeed, $N_1$ contains
pairs with arbitrary small negative values of $r$. Therefore, if
$\rho_0 < 0$ we can increase the right side of (7) as much as we
wish and get a contradiction with (7).

Clearly, $(0,I(x_0) \in S$. Thus, for each $z \leq 0$ we have
$(z,I(x_0)) \in S$. On the other hand $(0,I(x_0)) \in N$. Then for
each $z \leq 0$ by (7)

$$\rho_0 I(x_0) + (y_0^{*},z) \geq \rho_0 I(x_0)$$

Consequently, for all $z \leq 0$ we get $(y_0^{*},z) \geq 0$.
Therefore, $y_0^{*} \leq 0$.

For each $x \geq 0$ we have $(b - A x, I(x)) \in S$. Then from (7)
we get

\begin{equation}
\rho_0 I(x) + (y_0^{*},b-Ax) \geq \rho_0 I(x_0)
\end{equation}

for each $x \geq 0$.

Let us show that $\rho_0 > 0$. Indeed, if $\rho_0 = 0$, then from
(8) we get

\begin{equation}
(-y_0^{*},b-Ax) \leq 0
\end{equation}

for each $x \geq 0$.

Since (2) are strong simultaneous (9) contradicts Lemma 3.

Thus, $\rho_0 > 0$ and $y_0^{*} \leq 0$. Let $z_0^{*} = -
\frac{y_0^{*}}{\rho_0}$. Then $z_0^{*} \geq 0$ and from (8) we
have

\begin{equation}
I(x) + (z_0^{*},Ax - b) \geq I(x_0)
\end{equation}

for each $x \geq 0$.

If we put $x=x_0$ in (10) we get

$$(z_0^{*},Ax_0 - b) \geq 0$$

On the other hand $z_0^{*} \geq 0, Ax_0 \leq b$ and consequently
$(z_0^{*},Ax_0 - b) \leq 0$. Last two inequalities imply that

\begin{equation}
(z_0^{*},Ax - b) = 0
\end{equation}

Now (10) implies the second inequality in (6).

Let us prove the first inequality. Clearly, $(z^{*},Ax_0 - b) \leq
0$ for each $z^{*} \geq 0$. By using (11) we get

$$(z^{*},Ax_0 - b) \leq (z_0^{*},Ax_0 - b)$$

for each $z^{*} \geq 0$. Therefore,

$$I(x_0) + (z^{*}, Ax_0 - b) \leq I(x_0) + (z_0^{*}, Ax_0 - b)$$

for each $z^{*} \geq 0$. The first inequality of (6) is proved.

Now we state a theorem establishing the Kuhn - Tucker condition
for the problem (1),(2).\\

{\bf Theorem  2}. Suppose that the constraints (2) satisfy the
strong simultaneity condition. Then the necessary and sufficient
condition for the existence of a solution $x_0$ of the problem
(1),(2) is the existence of a nonzero linear functional $z^{*}_0
\geq 0$ such that the following conditions are held:

\begin{equation}
I'(x_0) + A^{*}z_0^{*} \geq 0
\end{equation}

\begin{equation}
(I'(x_0) + A^{*}z_0^{*},x_0) =  0
\end{equation}

\begin{equation}
Ax - b \geq 0, x_0 \geq  0
\end{equation}

\begin{equation}
(z_0^{*},Ax - b) = 0
\end{equation}

where $I'(x)$ is a gradient of $I(x)$, $A^{*}$ is the operator
adjoint to $A$.\\

{\bf Proof}. Due to Theorem 1, in order to prove theorem we have
to establish that the condition (6) is equivalent to the
conditions  (12)-(15).

Suppose that (6) is held. The second inequality of (6) means that
$x_0$ is a minimal point of convex functional $L(x,z_0^{*})$. By
the convex differentiability of a linear functional  for each $x
\geq 0$

$$(L'_{x}(x_0,z_0^{*}),x-x_0) \geq 0$$

Since  $L'_{x}(x_0,z_0^{*}) = I'(x_0) + A^{*}z_0^{*}$ we obtain
that for each $x \geq 0$

\begin{equation}
(I'(x_0) + A^{*}z_0^{*},x-x_0)\geq 0
\end{equation}

Consequently, $I'(x_0) + A^{*}z_0^{*} \geq 0$.

Put $x = 0$ in (16):

$$I'(x_0) + A^{*}z_0^{*},x_0)\leq 0 $$

On the other hand, $I'(x_0) + A^{*}z_0^{*},x_0)\geq 0$. Last two
inequalities imply (13).

First inequality of (6) implies that for each $z^{*} \geq 0$

\begin{equation}
 (z^{*},Ax_0 - b) \leq (z_0^{*},Ax_0 - b)
\end{equation}

and consequently, for each $z^{*} \geq 0$

\begin{equation}
 (z^{*},Ax_0 - b) \leq 0
\end{equation}

or $Ax_0 - b \leq 0$.

Now we get $(z_0^{*},Ax_0 - b) \geq 0$ by putting $z^{*}=0$ in
(17). On the other hand, $z_0^{*} \geq 0$, $Ax_0 - b \leq 0$ and
hence $(z_0^{*},Ax_0 - b) \leq 0$. Last two inequalities imply
(15).

Now suppose that (12)-(15) are held. From (12) we get that for all
$x \geq 0$

$$(I'(x_0) + A^{*}z_0^{*},x) \geq 0$$

Now by using (13) we get that for all $x \geq 0$

$$(I'(x_0) + A^{*}z_0^{*},x-x_0) \geq 0$$

In other words, for all $x \geq 0$

$$(L'_x(x_0,z_0^{*}),x-x_0) \geq 0$$

The last inequality is a necessary and sufficient condition for
$x_0$ to be a minimal point of $L(x,z_0^{*})$ for $x \geq 0$.
Therefore, for all $x \geq 0$ we get $L(x_0,z_0^{*}) \leq
L(x,z_0^{*})$. Thus, the right side of (6) is proved.

From (14) we get $(z^{*},Ax_0-b) \leq 0$ for all $z^{*} \geq 0$.
Now by (15), we get $(z^{*},Ax_0-b) \leq (z_0^{*},Ax_0-b)$ for all
$z^{*} \geq 0$. Therefore, $L(x_0,z^{*}) \leq L(x_0,z_0^{*})$ for
all $z^{*} \geq 0$. Thus, the left side of (6) also is proved.\\

{\bf Remark}. It can be readily shown that the strong simultaneity
condition (2) is  equivalent to the following condition

$$ 0 \in int(AK + b + P)$$

Clearly, $AK + b + P$ can have interior points even if $P$ has no
interior points. It means that the strong simultaneity condition
can be held in cases when Slater condition is not held.\\

{\bf Proposition}. In the case when $int P \not= \emptyset$, the
Slater
and the strong simultaneity conditions are equivalent.\\

{\bf Proof}. Suppose that the Slater condition is held: there is a
point $x_0 \geq 0$ such that $b - Ax_0 \in int P$. Then readily
the strong simultaneity condition is held.

Now let the strong simultaneity condition is held. Then there
exists a real number $\rho > 0$ such that for each $y \in
S_{\rho}$ ($S_{\rho}$ is a sphere with radius $\rho$ centered at
0) $b - Ax \geq y$, $x \geq 0$. Clearly, the strong simultaneity
condition can be written as

\begin{equation}
S_{\rho} \subset AK - b + P
\end{equation}

In order to prove that the Slater condition is held we show that
there exists a point $x_0 \geq 0$ and a real number $\rho_1 > 0$
such that

$$S_{\rho_1} \subset Ax_0 - b + P$$

It suffices to show that

$$int P \cap b - AK \not= \emptyset$$

Suppose the contrary: $int P \cap b - AK = \emptyset$. Since $P$
and $b - AK$ are convex, by separation theorem [2] there exists a
linear functional $z_0^{*} \in Y^{*}, z_0^{*} \not= 0$ such that

$$(z_0^{*},P) \leq 0 \leq z_0^{*},b - AK)$$

or equivalently, $(z_0^{*}, AK - b + P) \leq 0$. From (16) we get
$(z_0^{*},S_{\rho}) \leq 0$. Thus, $z^{*} = 0$. This is a
contradiction. The proof is completed.\\

\vspace{1cm}

{\bf References}\\

[1] Hurwitz L. , Udzava H. Stanford University Press (1958) \\

[2] Danford N., Schwartz J.T. Linear Operators, 1, (1988) \\

[3]  Burachik R.S., Jeyakumar V. Mathematical Programming, Springer,104,2-3 (2005). \\

[4] Evans L.C., Gomez D. Control, Optimization and Calculus of
Variations, v.8, (2002).\\

[5] Molanowski K. ,  Journal of  Applied mathematics and
Optimization, 25, 1, (1992).\\

[6] Chen S.Y., Wu S.Y., Journal of Computational and Applied
Mathematics, 213, 2, (2008).\\

[7] Ito S. Journal of Industrial and Management Optimization, v.6,
1, (2010).\\
\end{document}